\documentclass[12pt]{amsart}

\usepackage{amssymb}
\usepackage{latexsym}
\usepackage{amsmath}

\begin{document}

\def\z{{\bf Z}}
\def\n{{\bf N}}
\def\c{{\bf C}}
\def\s{{\bf S}}
\def\u{{\bf U}}
\def\vs{\bigskip}
\def\lvs{\vskip 6pt}

\thispagestyle{empty}
\null
\vskip 3truecm
{\normalsize \bf SUBFACTORS ASSOCIATED TO COMPACT KAC ALGEBRAS}
\vskip 1truecm
{\normalsize TEODOR BANICA}
\vskip 1truecm
\baselineskip=12pt
We construct inclusions of the form $(B_0\otimes P)^G\subset (B_1\otimes P)^G$, where $G$ is a compact quantum group of Kac type acting on an inclusion of finite dimensional $\c^*$-algebras $B_0\subset B_1$ and on a $II_1$ factor $P$. Under suitable assumptions on the actions of $G$, this is a subfactor, whose Jones tower and standard invariant can be computed by using techniques of A. Wassermann. The subfactors associated to subgroups of compact groups, to projective representations of compact groups, to finite quantum groups, to finitely generated discrete groups, to vertex models and to spin models are of this form.
\vskip 1truecm

\baselineskip=15pt

\newtheorem{defi}{Definition}[section]
\newtheorem{prop}{Proposition}[section]
\newtheorem{theo}{Theorem}[section]
\newtheorem{lemm}{Lemma}[section]

{\bf 1. INTRODUCTION}

\vs

There exist several constructions of subfactors using quantum groups and vice versa. The precise relation between this objects is far from being clear. In this paper we present a construction of subfactors using compact quantum groups of Kac type, which unifies six previously known constructions of subfactors (see the abstract). This is done by extending to quantum groups the following result.

\lvs

{\bf THEOREM 1.1.} {\it Let $G$ be a compact group. Let $P$ be a $II_1$ factor and let $G\to Aut(P)$ be a minimal action. Let $B_0\subset B_1$ be a Markov inclusion of finite dimensional von Neumann algebras and let $G\to Aut(B_1)$ be an action which leaves invariant $B_0$ and which is such that its restrictions to the centers of $B_0$ and $B_1$ are ergodic.

(i) $(B_0\otimes P)^G\subset (B_1\otimes P)^G$ is an inclusion of $II_1$ factors.

(ii) its Jones tower is $(B_1\otimes P)^G\subset  (B_2\otimes P)^G\subset (B_3\otimes P)^G\subset\cdots$, where $\{ B_i\}_{i\geq 1}$ are the algebras in the Jones tower for $B_0\subset B_1$, endowed with the canonical actions of $G$ coming from the action $G\to Aut(B_1)$.

(iii) its standard invariant is $\{ (B_i^\prime\cap B_j)^G\}_{0\leq i\leq j<\infty}$.}

\lvs

There are three main examples of such subfactors. The corresponding three main particular cases of theorem 1.1 can be found in \cite{w2} and are of course the source of inspiration for theorem 1.1. These are the Wassermann subfactors, the group-subgroup subfactors and the subfactors associated to crossed products by finite groups:
$$\left( P^G\subset (M_n \otimes P)^G\right)\,\simeq\, \left( ({\bf C}\otimes P)^G\subset (M_n\otimes P)^G\right)$$
$$\left( P^G\subset P^H\right)\,\simeq\, \left( ({\bf C}\otimes P)^G\subset (l^\infty(G/H)\otimes P)^G\right)$$
$$\left( P\subset P\rtimes G\right)\,\simeq\,  \left( (l^\infty (G)\otimes P)^G\subset ({\mathcal L} (l^2(G))\otimes P)^G\right)$$

We will extend theorem 1.1 to the compact quantum groups of Kac type. There is an obvious obstruction: the tensor product of two actions of a quantum group makes no sense in general. The point is that there exists a reasonable definition for the fixed point algebra of the ``non-existing'' (or maybe just ``bad'', see \S 3 below) tensor product and this is exactly what we need for stating the result (theorem 5.1 below).

Besides the above three classes of subfactors we may obtain also their ``quantum analogues''. For instance we may obtain Ocneanu's subfactor $P\subset P\rtimes G$, with $G$ a finite quantum group. Also, the diagonal subfactor associated to a finitely generated outer discrete group of automorphisms of a $II_1$ factor $\Gamma =<g_1,\ldots ,g_n>\subset Aut (Q)$ is isomorphic to a generalised Wassermann subfactor:
$$\left( \left\{ \begin{pmatrix}
g_1(q) & 0 & 0 \cr
0 & \cdots & 0 \cr
0 & 0 & g_n(q)
\end{pmatrix} \mid q\in Q\right\} \subset M_n(Q)  \right)\,\simeq\,  \left( (Q\rtimes\Gamma )^{\widehat{\Gamma}}\subset (M_n \otimes (Q\rtimes\Gamma ))^{\widehat{\Gamma}}\right)$$

Here $\widehat{\Gamma}$ acts on $Q\rtimes\Gamma$ via the dual of the action $\Gamma\subset Aut (Q)$ and on $M_n$ via the adjoint action of the representation $\oplus g_i :\widehat{\Gamma}\to {\bf C}^n$. At the combinatorial level, this was already pointed out in \cite{subf}. More generally, the generalised Wassermann subfactors provide spatial realisations for the Popa systems associated to representations of compact quantum groups \cite{subf}, under the Kac type assumption.

These are the four fundamental examples of subfactors of the form $(B_0\otimes P)^G\subset (B_1\otimes P)^G$. In \cite{sq} we prove that the subfactors associated to vertex models and to spin models (in the sense of V. Jones) are also of this form.

There are some natural questions about these subfactors. Let us us call ``analytical data'' the $II_1$ factor $P$ together with the minimal action $G\to Aut (P)$ and  ``algebraic data'' the Markov inclusion $B_0\subset B_1$ together with the action $G\to Aut(B_1)$ which leaves invariant $B_0$ and which is such that its restrictions to the centers of $B_0$ and $B_1$ are ergodic.

-- By \cite{u} any compact quantum group of Kac type acts minimally on at least one $II_1$-factor, so in particular it is part of some analytical data. Since the subfactor we construct is hyperfinite if and only if $P^G$ is hyperfinite, the remaining question about the analytical data is: does any compact quantum group of Kac type act minimally on a $II_1$ factor $P$ such that $P^G$ is hyperfinite? Note that this is weaker than the well-known question: does any discrete quantum group of Kac type act outerly on the hyperfinite $II_1$ factor?

-- Our result associates a Popa system to any algebraic data. It would be interesting to have more delicate results about this construction, like those in \cite{subf}. The first related question is: given an inclusion $B_0\subset B_1$, is it part of an algebraic data? This is true if $B_0=\c$: by \cite{aut} S. Wang's universal quantum group $G_{aut}(B_1)$ acts ergodically on $Z(B_1)$. In the general case one has to replace the quantum group $G_{aut}(B_1)$ from  \cite{aut} by a certain quantum group $G_{aut}(B_0\subset B_1)$ and the Temperley-Lieb algebra from \cite{aut} by the Fuss-Catalan algebra \cite{bj}. One gets that $B_0\subset B_1$ is part of an algebraic data if and only if it satisfies the necessary conditions that we find in \S 6 below. This will be done in a forthcoming paper.

-- In \cite{u2} Y. Ueda associates subfactors to $\s\u (2)_q$, which is a compact quantum group not of Kac type. It is not clear how to unify his construction of subfactors and ours.

The rest of the paper is organised as follows. \S 2 is a preliminary section. In \S 3 we define von Neumann algebras of the form $(B\otimes P)^G$ and we give a list of examples. In \S 4 we study the factoriality of such algebras. In \S 5 we prove the main result on inclusions of such factors and we give a list of examples (the semiduality of minimal actions, which is one of the main ingredients of the proof, is the subject of the appendix \S 7). In \S 6 we find an obstruction to non-integer values of the index.

I would like to thank M. Enock, S. Popa, Y. Ueda and A. Wassermann for useful discussions.

\vs

{\bf 2. COACTIONS AND COREPRESENTATIONS}
\vs

A compact quantum group of Kac type $G$ is an abstract object which may be described by three algebras: (1) the $*$-algebra ${\mathcal A}$ of ``representative functions on $G$''; this is a cosemisimple Hopf $*$-algebra whose square of the antipode is the identity, (2) the ${\bf C}^*$-algebra ${\bf A}$ of ``continuous functions on $G$''; this is a Woronowicz-Kac algebra, (3) the von Neumann algebra $A$ of ``$L^\infty$ functions on $G$''; this is a compact Kac algebra. See \cite{es}, \cite{wo}, \cite{bs} where more general objects (locally compact quantum groups of Kac type, compact matrix quantum groups, resp. locally compact quantum groups) are considered.

More precisely, we use the following formalism. Let ${\bf A}$ be a Woronowicz-Kac algebra with comultiplication $\delta :{\bf A}\to {\bf A}\otimes {\bf A}$ and antipode $\kappa :{\bf A}\to{\bf A}$. The canonical dense $*$-subalgebra ${\mathcal A}$ is a Hopf $*$-algebra, with the restrictions of $\delta$ and $\kappa$ as comultiplication and antipode and with counit denoted by $\varepsilon :{\mathcal A}\to {\bf C}$. Let $h:{\mathcal A}\to {\bf C}$ be the Haar trace and consider the left regular representation ${\mathcal A}\subset B(l^2({\mathcal A} ))$. The von Neumann algebra $A$ generated by the image of ${\mathcal A}$ is a compact Kac algebra, with the comultiplication, antipode and Haar trace given by the unique normal extensions of $\delta$, $\kappa$, $h$.

If $H$ is a finite dimensional Hilbert space and $u\in {\mathcal L}(H)\otimes A$ satisfies the coassociativity condition $(id\otimes\delta )u=u_{12}u_{13}$ and the smoothness condition $u\in {\mathcal L}(H)\otimes {\mathcal A}$ then the following four ``unitarity'' conditions are equivalent (cf. \cite{wo}):
$$uu^*=1,\,\,\, u^*u=1,\,\,\, u^t\overline{u}=1,\,\,\, \overline{u}u^t=1$$
If they are satisfied, $u$ is called a unitary corepresentation of $A$ on $H$.

All the von Neumann algebras in this paper are finite and their given faithful normal unital traces are usually denoted by $tr$. They are endowed with the scalar product $<x,y>=tr(y^*x)$. A coaction of $A$ on a finite von Neumann algebra $P$ is an injective morphism of von Neumann algebras $\pi :P\rightarrow P\otimes A$ satisfying the coassociativity condition $(\pi\otimes id)\pi =(id\otimes\delta )\pi$, the equivariance of the trace condition $(tr\otimes id)\pi =tr(.)1$ and the smoothness condition $\overline{{\mathcal P}}^w=P$, where ${\mathcal P}:=\pi^{-1}(P\otimes_{alg}{\mathcal A})$.

We use the following two constructions relating coactions and corepresentations.

\lvs

{\bf PROPOSITION 2.1.} {\it If $u\in{\mathcal L} (H)\otimes A$ is a unitary corepresentation and $\pi :P\to P\otimes A $ is a coaction then $\pi_u: x\mapsto u_{13}((id\otimes\pi )x)u_{13}^*$ defines a coaction of $A$ on ${\mathcal L}(H)\otimes P$.}

\lvs

{\bf PROOF.} The coassociativity (resp. smoothness) of $\pi_u$ is clear from the coassociativity (resp. smoothness) of $\pi$ and of $u$. If $\{ e_i\}$ is a basis of $H$ then
$$(tr\otimes tr\otimes id)(\pi_u (e_{ij}\otimes p))=(tr\otimes tr\otimes id)(\sum_{ab}e_{ab}\otimes ((1\otimes u_{ai})\pi (p) (1\otimes u^*_{bj})))=$$
$$n^{-1}\sum_a u_{ai}((tr\otimes id)\pi (p))u^*_{aj}=n^{-1}tr(p)(u^t\overline{u} )_{ij}=(tr\otimes tr)(e_{ij}\otimes p)1$$
for any $i,j,p$. Thus the trace is $\pi_u$-equivariant. 

\lvs

{\bf PROPOSITION 2.2.} {\it If $\beta :B\to B\otimes A$ is a coaction on a finite dimensional algebra then its image $u_\beta$ by the canonical linear isomorphism ${\mathcal L}(B,B\otimes A)\simeq {\mathcal L}(B)\otimes A$ is a unitary corepresentation. If $\{ b_i\}$ is an orthonormal basis of $B$ and we write $\beta (b_i)=\sum_j b_j\otimes u_{ji}$ then $u_\beta =\sum_{ij}e_{ij}\otimes u_{ij}$.}

\lvs

{\bf PROOF.} The coassociativity (resp. smoothness) of $u_\beta$ is clear from the coassociativity (resp. smoothness) of $\beta$. We have
$$\sum_j u_{ji}u_{jk}^*=(tr\otimes id)(\sum_{js} b_jb_s^*\otimes u_{ji}u_{sk}^*)=(tr\otimes id)\beta (b_ib_k^*)=tr(b_ib_k^*)1=\delta_{i,k}$$
for any $i,k$ and it follows that $u^t\overline{u}=1$. 

\vs

{\bf 3. FIXED POINT ALGEBRAS}

\vs

If $\beta :B\to B\otimes A$ and $\pi :P\to P\otimes A$ are two coactions there is no way of constructing a coaction $B\otimes P\to B\otimes P\otimes A$ which reasonably corresponds to the notion of ``tensor product of $\beta$ and $\pi$''. More dramatically, one cannot construct a subalgebra of $B\otimes P$ which reasonably corresponds to the notion of ``fixed point algebra under the tensor product of $\beta$ and $\pi$''. However, we will see in this section that there exists one reasonable construction, which associates to any anticoaction $\beta :B\to B\otimes A$ and to any coaction $\pi :P\to P\otimes A$ a certain subalgebra $(B\otimes P)^{\beta\otimes\pi}$ of $B\otimes P$.

\lvs

{\bf DEFINITION 3.1.} {\it An anticoaction of $A$ on a finite dimensional von Neumann algebra $B$ is a map $\beta :B\to B\otimes A$ satisfying one of the following six equivalent conditions:

(i) the map $\!\!{\ }^o\beta :=(o\otimes id)\beta\, o:B^o\rightarrow B^o\otimes A$ is a coaction, where $B^{o}$ is the opposite algebra and $o:B\leftrightarrow B^{o}$ are the canonical maps.

(ii)  the map $\!\!{\ }^t\beta :=(t\otimes id)\beta\, t:B\rightarrow B\otimes A$ is a coaction, where $t:B\to B$ is a transposition.

(iii) the map $\beta^\kappa =(id\otimes\kappa )\beta :B\rightarrow B\otimes A$ is a coaction of the Kac algebra $(A,\sigma\delta ,\kappa ,h)$, where $\sigma$ is the flip.

($i^\prime$) (resp. ($ii^\prime$), resp. ($iii^\prime$)) the map $\beta$ satisfies the coassociativity, equivariance of the trace and smoothness assumptions and the map $\!\!{\ }^o\beta$ (resp. $\!\!{\ }^t\beta$, resp. $\beta^\kappa$) is an injective $*$-morphism.}

\lvs

Here the three equivalences $(j)\Leftrightarrow (j^\prime )$ are clear from definitions and the equivalence $(i^\prime )\Leftrightarrow (ii^\prime )\Leftrightarrow (iii^\prime )$ is clear from the fact that $o$, $t$, $\kappa$ are antimorphisms. 

If $A=L^\infty (G)$ with $G$ a compact group then coactions and anticoactions coincide (cf. $(iii^\prime )$ and the fact that the antipode is a $*$-morphism). When $A$ is non-commutative, this is no longer true.

The canonical linear map $o:B\to B^o$ induces an isomorphism ${\mathcal L}(B)\simeq {\mathcal L}(B^o)$. Thus proposition 2.2 may be used for associating unitary corepresentations on $B$ to anticoactions on $B$. More precisely, if $\beta :B\to B\otimes A$ is an anticoaction,  $\{ b_i\}$ is an orthonormal basis of $B$ and we write $\beta (b_i)=\sum_j b_j\otimes u_{ji}$ then $u_\beta =\sum_{ij}e_{ij}\otimes u_{ij}$ is a unitary corepresentation.

\lvs

{\bf DEFINITION 3.2.} {\it If $\pi :P\to P\otimes A$ is a coaction on a von Neumann algebra and $\beta :B\to B\otimes A$ is an anticoaction on a finite dimensional von Neumann algebra we define a coassociative linear map $\beta\otimes\pi :B\otimes P\to B\otimes P\otimes A$ by $b\otimes p\mapsto \pi (p)_{23}\beta (b)_{13}$.}

\lvs

If $Q$ is a von Neumann algebra and $\Gamma :Q\to Q\otimes A$ is a coassociative linear map we define a linear map $E_\Gamma :=(id\otimes h)\Gamma :Q\to Q$. Then $Im(E_\Gamma)=Q^\Gamma$, the fixed point algebra of $Q$ under $\Gamma$ (the inclusion $Im(E_\Gamma )\subset Q^\Gamma$ follows from the computation
$$\Gamma E_\Gamma =(id\otimes id\otimes h)(\Gamma\otimes id )\Gamma =(id\otimes id\otimes h)(id\otimes\delta )\Gamma =(id\otimes h(.)1)\Gamma =E_\Gamma \otimes 1$$
and the reverse inclusion is clear). If $\Gamma$ is a coaction or an anticoaction then $Q^\Gamma$ is a von Neumann algebra and $E_\Gamma$ is the conditional expectation onto it.
\lvs

{\bf LEMMA 3.1.} {\it We have $E_{\pi_{u_\beta}}(\lambda\otimes id)=(\lambda\otimes id)E_{\beta\otimes\pi}$, where $\lambda :B\to {\mathcal L}(B)$ is the left regular representation $\lambda (x):y\mapsto xy$.}

\lvs

{\bf PROOF.} Let $\{ b_i\}$ be an orthonormal basis of $B$ and write $\beta (b_i)=\sum_j b_j\otimes u_{ji}$. We will use many times the formula $b=\sum_xb_xtr(bb_x^*)$ for any $b\in B$. By linearity it is enough to prove the formula on $b\otimes p$ with $b,p$ arbitrary, i.e. to prove that
$$(id\otimes id\otimes h)((u_\beta )_{13}(\lambda (b)\otimes \pi (p))(u^*_\beta )_{13})= \sum_{sx} tr(bb_x^*)\lambda (b_s)\otimes (id\otimes h)(\pi (p)(1\otimes u_{s x}))$$
By applying $id\otimes f$ with $f\in P^*$ arbitrary, we want to prove that
$$(id\otimes h)(u_\beta (\lambda (b)\otimes \xi )u^*_\beta )=\sum_{sx} tr(bb_x^*)\lambda (b_s ) h(\xi u_{s x})$$
where $\xi =(f\otimes id)\pi (p)$. The left term is $\sum_{ijkl}e_{ij}\lambda (b_x)e_{kl}h(\xi u^*_{lk}u_{ij})$, so the above formula is obtained from the one below by applying $id\otimes h(\xi .)$
$$\sum_{ijkl}e_{ij}\lambda (b)e_{kl}\otimes u^*_{lk}u_{ij}=\sum_{sx} tr(bb_x^*) \lambda (b_s ) \otimes u_{s x}$$
By applying $<.b_l,b_i>\otimes id$ with $i,l$ arbitrary we want to prove that
$$\sum_{jk}tr(bb_kb_j^*)u^*_{lk}u_{ij}=\sum_{sx} tr(bb_x^*) tr(b_i^*b_s b_l)u_{sx}$$
The above formula is obtained from the one below by applying $tr(b.)\otimes id$
$$\sum_{jk}b_kb_j^*\otimes u^*_{lk}u_{ij}=\sum_{sx} tr(b_i^*b_s b_l)b_x^*\otimes u_{sx}$$
The right term is $\beta^\kappa (\sum_s tr(b_sb_lb_i^*)b_s^*))=\beta^\kappa (b_lb_i^*)$ and the left term is $\beta^\kappa(b_l)\beta^\kappa (b_i)^*$, where $\beta^\kappa :b_i\mapsto\sum_j b_j\otimes u_{ij}^*$ is the $*$-morphism in definition 3.1 (iii).

\lvs

{\bf THEOREM 3.1.} {\it $(B\otimes P)^{\beta\otimes\pi}$ is a von Neumann subalgebra of $B\otimes P$ and $E_{\beta\otimes\pi}$ is the conditional expectation onto it.}

\lvs

{\bf PROOF.} If we regard $\lambda\otimes id$ as an inclusion then we have the following equalities between subsets of ${\mathcal L}(B)\otimes P$:
$$(B\otimes P)^{\beta\otimes\pi}=E_{\beta\otimes\pi}(B\otimes P)=E_{\pi_{u_\beta}}(B\otimes P)=(B\otimes P)\cap ({\mathcal L}(B)\otimes P)^{\pi_{u_\beta}}$$
The first equality is clear, the second one follows from lemma 3.1 and the third one follows from the fact that $E_{\pi_{u_\beta}}$ is a conditional expectation onto its image $({\mathcal L}(B)\otimes P)^{\pi_{u_\beta}}$. Since both $B\otimes P$ and $({\mathcal L}(B)\otimes P)^{\pi_{u_\beta}}$ are von Neumann algebras, it follows that their intersection $(B\otimes P)^{\beta\otimes\pi}$ is a von Neumann algebra. Also from the commutativity of the above diagram we get that $E_{\beta\otimes\pi}$ is the restriction of $E_{\pi_{u_\beta}}$, so it is a conditional expectation onto its image.

\lvs

{\bf EXAMPLE 3.1.} Let $G$ be a compact group and let $p:G\to Aut(P)$ and $b:G\to Aut(B)$ be two actions of it. Then the corresponding maps $\pi :P\to P\otimes L^\infty (G)$ and $\beta :B\to B\otimes L^\infty (G)$ are both coactions and anticoactions and we have the following equality between subalgebras of $B\otimes P$:
$$(B\otimes P)^{\beta\otimes\pi} =(B\otimes P)^{b\otimes p}$$
Indeed, it is easy to see that $\beta\otimes\pi$ is the coaction of $L^\infty (G)$ associated to $b\otimes p$.
\lvs

{\bf EXAMPLE 3.2.} Let $A$ be a finite dimensional Kac algebra. Since $\kappa$ is a transposition of $A$, the map $\!\!{\ }^\kappa\delta =(\kappa\otimes id )\delta\kappa$ in definition 3.1 (ii) is an anticoaction. If $\pi :P\to P\otimes A$ is a coaction then the following equality between subalgebras of $A\otimes P$ holds
$$(A\otimes P)^{\!\!{\ }^\kappa\delta\otimes\pi}=\sigma\pi (P)$$
where $\sigma :P\otimes A\to A\otimes P$ is the flip. For, let $\{ u_{ij}\}$ be an orthonormal basis of $A$ consisting of coefficients of irreducible corepresentations (see e.g. \cite{wo}). For any $p\in P$ we use the notation $\pi (p)=\sum_{uij}p^u_{ij}\otimes u_{ij}$. From the coassociativity of $\pi$ we get that $\pi (p_{ij}^u)=\sum_k p_{kj}^u\otimes u_{ki}$ for any $u,i,j$. We get succesively that
$$\!\!{\ }^\kappa\delta (u_{ij})=(\kappa\otimes id )\delta (u_{ji}^*)=\sum_k \kappa (u_{jk}^*)\otimes u_{ki}^* =\sum_k u_{kj}\otimes u_{ki}^*$$
$$(\!\!{\ }^\kappa\delta\otimes \pi )(u_{ij}\otimes p)=\pi (p)_{23}\!\!{\ }^\kappa\delta (u_{ij})_{13}=\sum_{vskl}u_{kj}\otimes p_{sl}^v\otimes v_{sl}u_{ki}^*$$
$$E_{\!\!{\ }^\kappa\delta\otimes \pi}(u_{ij}\otimes p)=\sum_k u_{kj}\otimes p_{ki}^u=\sigma\pi (p_{ij}^u)$$
Thus $(A\otimes P)^{\!\!{\ }^\kappa\delta\otimes\pi}\subset \sigma\pi (P)$. The other inclusion follows from the following formula:
$$E_{\!\!{\ }^\kappa\delta\otimes \pi}(\sigma\pi (p))=E_{\!\!{\ }^\kappa\delta\otimes \pi}(\sum_{uij}u_{ij}\otimes p^u_{ij})=\sum_{uijkl}u_{kj}\otimes p_{lj}^uh(u_{li}u_{ki}^*)=\sigma\pi (p)$$
\lvs

{\bf EXAMPLE 3.3.} Denote by $\iota :{\bf C}\to {\bf C}\otimes A$ the trivial coaction. If $v\in M_n(A)$ is a corepresentation then the map $\iota_v :M_n\to M_n\otimes A$ given by proposition 2.1 is a coaction, hence the map $\!\!{\ }^t\iota_v :M_n\to M_n\otimes A$ given by definition 3.1 (ii) is an anticoaction. If $\pi :P\to P\otimes A$ is a coaction then the following subalgebras of $M_n\otimes P$ are equal
$$(M_n\otimes P)^{\!\!{\ }^t\iota_v\otimes\pi}=(M_n\otimes P)^{\pi_{\overline{v}}}$$
Indeed, the maps $\pi_{\overline{v}}$ and $\!\!{\ }^t\iota_v\otimes\pi$ are given by
$$\pi_{\overline{v}}:e_{ij}\otimes p\mapsto {\overline{v}}_{13}(e_{ij}\otimes \pi (p))v_{13}^t=\sum_{ab}e_{ab}\otimes ((1\otimes v^*_{ai})\pi (p)(1\otimes v_{bj}))$$
$$\!\!{\ }^t\iota_v\otimes\pi :e_{ij}\otimes p\mapsto \pi (p)_{23}((t\otimes id)(v(e_{ji}\otimes 1)v^*))_{13}=\sum_{ab}e_{ab}\otimes (\pi (p)(1\otimes v_{bj}v^*_{ai}))$$
and by using the trace property of $h$ we get from this that $E_{\pi_{\overline{v}}}=E_{\!\!{\ }^t\iota_v\otimes\pi}$.
\lvs

{\bf EXAMPLE 3.4.} Let $\Gamma$ be a discrete group and consider the compact Kac algebra $L(\Gamma )$. If $g_1,\ldots ,g_n$ are elements of $\Gamma$ then $v=diag (u_{g_i})$ is a corepresentation of it, where $g\mapsto u_g$ the canonical map $\Gamma\to L(\Gamma )$. Its associated anticoaction $\!\!{\ }^t\iota_v:M_n\to M_n\otimes L(\Gamma )$ is given by $e_{ij}\mapsto e_{ij}\otimes u_{g_jg_i}^{-1}$ (see example 3.3). Let $G\to Aut(Q)$ be an action on a von Neumann algebra and consider the dual coaction $\pi :P\to P\otimes L(\Gamma )$ on the crossed product $P=Q\rtimes\Gamma$. Then the following subalgebras of $M_n\otimes P$ are equal:
$$(M_n\otimes P)^{\!\!{\ }^t\iota_v\otimes\pi} =v(M_n\otimes Q)v^*$$
where ``$Q$'' denotes the canonical copy of $Q$ inside $P$ and ``$v$'' denotes the image of $v\in M_n\otimes L(\Gamma )$ in $M_n\otimes P$ by the canonical inclusion. Indeed, if $p=\sum_{ij}e_{ij}\otimes p_{ij}$ is an arbitrary element in $M_n\otimes P$ then $(\!\!{\ }^t\iota_v\otimes\pi )(p)=\sum_{ij}e_{ij}\otimes (\pi (p_{ij})(1\otimes u_{g_jg_i^{-1}}))$, so $p$ is in $(M_n\otimes P)^{\!\!{\ }^t\iota_v\otimes\pi}$ if and only if $\pi (p_{ij})=p_{ij}\otimes u_{g_ig_j^{-1}}$ for any $i,j$. Since $\pi (u_h)=u_h\otimes u_h$ for any $h\in\Gamma$, this is equivalent to $u_{g_i^{-1}}p_{ij}u_{g_j}\in P^\pi =Q$ for any $i,j$, hence to $p_{ij}\in u_{g_i}Qu_{g_j^{-1}}$ for any $i,j$.

\vs

{\bf 4. FIXED POINT FACTORS}

\vs

Let $\pi :P\rightarrow P\otimes A$ be a coaction. If $u\in{\mathcal L} (H)\otimes A$ is a finite dimensional unitary corepresentation, an eigenmatrix for $u$ is an element $M\in {\mathcal L} (H)\otimes P$ satisfying $(id\otimes\pi )M=M_{12}u_{13}$ in ${\mathcal L} (H)\otimes P\otimes  A$. A coaction $\pi$ is called semidual if each corepresentation has a unitary eigenmatrix. The canonical coaction $\delta$ on $A$ is clearly semidual -- each corepresentation is a unitary eigenmatrix for itself. It follows that any coaction which contains equivariantly a copy of $A$ (such a coaction is said to be dual) is semidual.

\lvs

{\bf LEMMA 4.1.} {\it If $\pi :P\to P\otimes A $ is a semidual coaction and $\beta :B\to B\otimes A $ is an anticoaction then
$$\begin{matrix}
P &\subset & B\otimes P\cr 
\cup &\ &\cup \cr 
P^\pi &\subset & (B\otimes P)^{\beta\otimes\pi}\cr
\end{matrix}$$
is non-degenerate commuting square.}

\lvs

{\bf PROOF.} The commuting square condition is clear from the formulae of the vertical conditional expectations. Let $\{ b_i\}$ be an orthonormal basis of $B$, write $\beta (b_i)=\sum_j b_j\otimes u_{ji}$ and consider the corresponding corepresentation $u_\beta =\sum_{ij}e_{ij}\otimes u_{ij}$. Let $m$ be a unitary eigenmatrix for $u_\beta$. Consider for any $i$ the element $A_i=\sum b_k^*\otimes m_{ik}$. Then
$$(\beta\otimes\pi )(A_i)=\sum_{kst}(m_{is}\otimes u_{sk})_{23}(b_t^*\otimes u_{tk}^*)_{13}=\sum_{kst}b_t\otimes m_{is}\otimes u_{sk}u_{tk}^*=A_i\otimes 1$$
so $A_i\in (B\otimes P)^{\beta\otimes\pi}$ for any $i$. For any $s$ and any $p\in P$ we have
$$\sum_i A_i^*(1\otimes m_{is}p)=\sum_{ik}b_k\otimes m_{ik}^*m_{is}p=b_s\otimes p$$
so $b_s\otimes p\in sp\{ (B\otimes P)^{\beta\otimes\pi} (1\otimes P)\}$ and this proves the non-degeneracy. 

\lvs

A coaction $\pi :P\to P\otimes A$ is said to be minimal if $(P^\pi )^\prime\cap P={\bf C}$ and if it faithful in the following sense: ${\overline{sp}^w\,} \{ (\phi\otimes id)\pi (p)\mid \phi\in P_*,\, p\in P\} = A$.

\lvs

{\bf THEOREM 4.1.} {\it If $\pi :P\to P\otimes A$ is a minimal coaction and $\beta :B\to B\otimes A $ is an anticoaction then $(B\otimes P)^{\beta\otimes\pi}$ is a factor if and only if $Z(B)\cap B^\beta ={\bf C}$.}

\lvs

{\bf PROOF.} We prove that the following subalgebras of $B\otimes P$ are equal:
$$Z((B\otimes P)^{\beta\otimes\pi}) = (Z(B)\cap B^\beta)\otimes 1$$
Let $x$ be in the algebra on the left. Then $x$ has to commute with $1\otimes P^\pi$, so by minimality it has to be of the form $b\otimes 1$. Thus $x$ commutes with $1\otimes P$. But $x$ commutes by definition with $(B\otimes P)^{\beta\otimes\pi}$ and from the non-degeneracy of the square in lemma 4.1 we get that $x$ commutes with $B\otimes P$ and in particular with $B\otimes 1$. Thus $b\in Z(B)\cap B^\beta$. The other inclusion is obvious. 

\vs

{\bf 5. FIXED POINT SUBFACTORS}

\vs

We will use the following known lemma ((i) is from \cite{p3} and (ii) is from \cite{w2}). 

\lvs

{\bf LEMMA 5.1.} {\it Consider two commuting squares of finite von Neumann algebras
$$\begin{matrix}
F &\subset & E &\subset & D\cr 
\cup &\ &\cup &\ &\cup\cr 
A &\subset & B  &\subset & C \cr
\end{matrix}$$

(i) if the square on the left and the big square are non-degenerate, so is the square on the right.

(ii) if both squares are non-degenerate, $F\subset E\subset D$ is the basic construction and its Jones projection $e\in D$ belongs to $C$ then the square on the right is the basic construction for the square on the left.}

\lvs

{\bf PROOF.} (i) is clear from $D={\overline{sp}^w\,} CF={\overline{sp}^w\,}CBF={\overline{sp}^w\,} CE$. For (ii), let $\Psi :D\to C$ be the expectation. We have that $E={\overline{sp}^w\,} FB={\overline{sp}^w\,} BF$ by non-degeneracy and $D={\overline{sp}^w\,} EeE$ by the basic construction so we get that
$$C=\Psi (D)=\Psi ({\overline{sp}^w\,} EeE)=\Psi ({\overline{sp}^w\,} BFeFB)=$$
$$\Psi ({\overline{sp}^w\,} BeFB)={\overline{sp}^w\,} Be \Psi (F)B={\overline{sp}^w\,} Be AB={\overline{sp}^w\,} BeB$$
Thus $C$ is generated by $B$ and $e$. 

\lvs

{\bf LEMMA 5.2.} {\it If $\beta :B\to B\otimes A$ is a coaction then
$$\begin{matrix}
A &\subset & B\otimes A\cr 
\cup &\ &\ \uparrow\beta\cr 
{\bf C} &\subset & B\cr
\end{matrix}$$
is a non-degenerate commuting square.}

\lvs

{\bf PROOF.} From the $\beta$-equivariance of the trace we get that the inclusion on the left commutes with the traces, i.e. that this is a commuting diagram of finite von Neumann algebras. From the formula of the expectation $E_{\beta}=(id\otimes h)\beta$ we get that this is a commuting square. Choose an orthonormal basis $\{ b_i\}$ of $B$, write $\beta :b_i\mapsto\sum_jb_j\otimes u_{ji}$ and consider the corresponding unitary corepresentation $u_\beta=\sum e_{ij}\otimes u_{ij}$. Then for any $k$ and any $a\in A$ we have
$$\sum_i\beta (b_i)(1\otimes u_{ki}^*a)=\sum_{ij}b_j\otimes u_{ji}u_{ki}^*a=\sum_{ij}b_j\otimes \delta_{j,k}a=b_k\otimes a$$
so $b_k\otimes a\in sp\,\beta (B)(1\otimes A)$. It follows that $B\otimes A=sp\, \beta (B)(1\otimes A)$, i.e. that our commuting square is non-degenerate. 

\lvs

{\bf PROPOSITION 5.1.} {\it Let $B_0\subset B_1$ be a Markov inclusion of finite dimensional von Neumann algebras and let
$$B_0\subset B_1\subset_{e_1} B_2=<B_1,e_1> \subset_{e_2} B_3=<B_2,e_2> \subset_{e_3} \cdots$$
be its Jones tower. If $\beta_1:B_1\to B_1\otimes A$ is a coaction (resp. anticoaction) leaving $B_0$ invariant then there exists a unique sequence $\{\beta_i\}_{i\geq 0}$ of coactions (resp. anticoactions) $\beta_i:B_i\to B_i\otimes A$ such that each $\beta_i$ extends $\beta_{i-1}$ and leaves invariant the Jones projection $e_{i-1}$.}

\lvs

{\bf PROOF.} By taking opposite inclusions we see that the assertion for anticoactions is equivalent to the one for coactions, that we will prove now. The unicity part is clear from $B_i=<B_{i-1},e_{i-1}>$. For the existence, by applying lemma 5.1 (i) to the diagram
$$\begin{matrix}
A &\subset & B_0\otimes A &\subset & B_1\otimes A\cr 
\cup &\ &\ \uparrow\beta_0 &\ &\ \uparrow\beta_1\cr 
{\bf C} &\subset & B_0 &\subset & B_1\cr
\end{matrix}$$
we get that the square on the right is a non-degenerate. By performing basic constructions to it we get a sequence
$$\begin{matrix}
B_0\otimes A &\subset & B_1\otimes A &\subset & B_2\otimes A &\subset & B_3\otimes A & \subset &\cdots \cr 
\uparrow\beta_0 &\ &\ \uparrow\beta_1 &\ &\ \uparrow\beta_2 &\ &\ \uparrow\beta_3\cr 
B_0 &\subset & B_1 &\subset & B_2 &\subset & B_3 & \subset &\cdots \cr
\end{matrix}$$
It is easy to see from definitions that the $\beta_i$'s are coactions, that they extend each other and that they leave invariant the Jones projections. 

\lvs

{\bf THEOREM 5.1.} {\it Let $\pi :P\to P\otimes A $ be a minimal coaction on a $II_1$ factor $P$, let $B_0\subset B_1$ be a Markov inclusion of finite dimensional von Neumann algebras and let $\beta_1 :B_1\to B_1\otimes A $ be a anticoaction which leaves $B_0$ invariant and which is such that $Z(B_0)\cap B_0^{\beta_0}=Z(B_1)\cap B_1^{\beta_1}={\bf C}$, where $\beta_0 :B_0\to B_0\otimes A$ is the restriction of $\beta_1$.

(i) $(B_0\otimes P)^{\beta_0\otimes\pi}\subset  (B_1\otimes P)^{\beta_1\otimes\pi}$ is an inclusion of $II_1$ factors.

(ii) its Jones tower is $(B_1\otimes P)^{\beta_1\otimes\pi}\subset  (B_2\otimes P)^{\beta_2\otimes\pi}\subset  (B_3\otimes P)^{\beta_3\otimes\pi}\subset\cdots$, where the $\beta_i$'s are the canonical anticoactions on the algebras $\{B_i\}_{i\geq 1}$ in the Jones tower for $B_0\subset B_1$.

(iii) its lattice of higher relative commutants is
$$\begin{matrix}
{\bf C} &\subset & B_0^\prime\cap B_1^{\beta_1} & \subset & B_0^\prime\cap B_2^{\beta_2}& \subset & B_0^\prime\cap B_3^{\beta_3}& \subset\ \cdots\cr
&\ &\ \cup &\ &\cup &\ &\cup \cr
&\ &\ {\bf C} & \subset & B_1^\prime\cap B_2^{\beta_2}& \subset & B_1^\prime\cap B_3^{\beta_3}& \subset\ \cdots\cr
&\ &\ &\ &\cup &\ &\cup \cr
&\ &\ \ & \ & {\bf C} & \subset & B_2^\prime\cap B_3^{\beta_3}& \subset\ \cdots\cr
&\ &\ &\ &\ &\ &\cup \cr
&\ &\ &\ &\ &\ &\cdots &\ \ \ \ \ \cdots\cr
\end{matrix}$$}

\lvs

{\bf PROOF.} (i) follows from theorem 4.1. Consider the following diagram
$$\begin{matrix}
P &\subset & B_i\otimes P &\subset & B_j\otimes P\cr 
\cup &\ &\cup &\ &\cup \cr 
P^\pi &\subset & (B_i\otimes P)^{\beta_i\otimes\pi} &\subset & (B_j\otimes P)^{\beta_j\otimes\pi}\cr
\end{matrix}$$
By lemma 4.1 the big square and the square on the left are non-degenerate commuting squares. Thus lemma 5.1 (i) applies and shows that the square on the right is a non-degenerate commuting square. Consider now the following sequence of non-degenerate commuting squares
$$\begin{matrix}
B_0\otimes P &\subset & B_1\otimes P &\subset &
B_2\otimes P &\subset & \cdots\cr 
\cup &\ &\cup &\ &\cup &\ &\ \cr 
(B_0\otimes P)^{\beta_0\otimes\pi} &\subset & (B_1\otimes P)^{\beta_1\otimes\pi} &\subset & (B_2\otimes P)^{\beta_2\otimes\pi} &\subset & \cdots\cr
\end{matrix}$$
Since the Jones projections live in the lower line, lemma 5.1 (ii) applies and shows that this is a sequence of basic constructions for non-degenerate commuting squares. In particular the lower line is a sequence of basic constructions and we get (ii). It remains to prove (iii), i.e. that for $0\leq i\leq j$ the following subalgebras of $B_j\otimes P$ are equal
$$((B_i\otimes P)^{\beta_i\otimes\pi})^\prime \cap (B_j\otimes P)^{\beta_j\otimes\pi} = (B_i^\prime\cap B_j^{\beta_j})\otimes 1$$
Let $x$ be in the algebra on the left. Then $x$ has to commute with $1\otimes P^\pi$, so by minimality it has to be of the form $b\otimes 1$. Thus $x$ commutes with $1\otimes P$. But $x$ commutes by definition with $(B_i\otimes P)^{\beta_i\otimes\pi}$ and from the non-degeneracy of the square on the left in the above diagram consisting of two squares we get that $x$ commutes with $B_i\otimes P$ and in particular with $B_i\otimes 1$. Thus $b\in B_i^\prime\cap B_j$, so $x=b\otimes 1$ is in the algebra on the right. The other inclusion is obvious. 

\lvs

{\bf EXAMPLE 5.1.} For $A=L^\infty (G)$ with $G$ a compact group we get exactly the theorem 1.1 (cf. example 3.1).
\lvs

{\bf EXAMPLE 5.2.} If $A$ is a finite dimensional Kac algebra and $\pi$ is a minimal coaction of it on a $II_1$ factor $P$ then the subfactor $P^\pi\subset P$ is isomorphic to the subfactor $P^\pi\subset (A\otimes P)^{\!\!{\ }^\kappa\delta\otimes\pi}$ (cf. example 3.2). By applying to it theorem 3.1 we get the well-known result on the Jones towers and standard invariants of such subfactors.
\lvs

{\bf EXAMPLE 5.3.} If $v\in M_n(A)$ is a corepresentation and $\pi :P\to P\otimes A$ is a minimal coaction on a $II_1$ factor $P$ then the subfactor $P^\pi\subset (M_n\otimes P)^{\pi_v}$ is isomorphic to the subfactor $P^\pi\subset (M_n\otimes P)^{\!\!{\ }^t\iota_{\overline{v}}\otimes\pi}$ (cf. example 3.3). By applying to it theorem 5.1 we get after an easy computation that its standard invariant is the following lattice $R(v)$:
$$\begin{matrix}
{\bf C} &\subset & End(v) & \subset & End(\overline{v}\otimes v)& \subset & End(v\otimes\overline{v}\otimes v)& \subset\ \cdots\cr
&\ &\ \cup &\ &\cup &\ &\cup \cr
&\ &\ {\bf C} & \subset & End(\overline{v})& \subset & End(v\otimes\overline{v})& \subset\ \cdots\cr
&\ &\ &\ &\cup &\ &\cup \cr
&\ &\ \ & \ & {\bf C} & \subset & End(v)& \subset\ \cdots\cr
&\ &\ &\ &\ &\ &\cup \cr
&\ &\ &\ &\ &\ &\cdots &\ \ \ \ \ \cdots\cr
\end{matrix}$$
The relation with \cite{subf} is as follows. If $({\bf A},v)$ is a finitely generated Woronowicz-Kac algebra then $({\bf A},v^*)$ is also a finitely generated Woronowicz-Kac algebra (with comultiplication $\sigma\delta$ and antipode $\kappa$) and it is easy to see that the lattice $R(v^*)$ is isomorphic to the lattice $L(v)$ as defined in \cite{subf}.
\lvs

{\bf EXAMPLE 5.4.} Let $\Gamma =<g_1,...,g_n>\subset Aut(Q)$ be a finitely generated outer discrete group of automorphisms of a $II_1$ factor $Q$. With the notations in example 3.4 we get that the generalised Wassermann subfactor $P^\pi\subset (M_n\otimes P)^{\pi_v}$ is isomorphic to the inclusion
$$Q\to M_n(Q),\,\,\, q\mapsto v(1\otimes q)v^*=diag(u_{g_i}qu_{g_i}^*)=diag(g_i(q))$$
which is the diagonal subfactor associated to $\Gamma$. By applying to it theorem 5.1 we get the well-known result on the Jones towers and standard invariants of such subfactors.

\vs

{\bf 6. VALUES OF THE INDEX}

\vs

We recall from \cite{aut} that given a ${\bf C}^*$-algebra $B$ of dimension $n<\infty$, for a faithful unital trace $tr:B\to {\bf C}$ the following conditions are equivalent:

-- $tr$ is the restriction of the unique trace of ${\mathcal L} (B)$, via the embedding $B\subset {\mathcal L}  (B)$ given by the left regular representation.

-- ${\bf C}\subset B$ is a Markov inclusion.

-- if $\phi :B\simeq \bigoplus_{\gamma =1}^sM_{m_\gamma}$ is a decomposition of $B$ as a multimatrix algebra, then the weights $\lambda_\gamma :=tr(\phi^{-1} (Id_{M_{m_{\gamma}}}))$ of $tr$ are given by $\lambda_\gamma =n^{-1}m_\gamma^2$ for any $\gamma$.

-- $\mu\mu^*=n\cdot id$, where $\mu :B\otimes B\to B$ is the multiplication and where the adjoint $\mu^*:B\to B\otimes B$ is taken with respect to the Hilbert space structure coming from $tr$.

This distinguished trace is called the canonical trace of $B$.

\lvs

{\bf LEMMA 6.1} {\it Let $(B,tr)$ be a finite dimensional finite von Neumann algebra and consider the element $\xi =\mu\mu^* (1_B)\in B$.

(i) For any anticoaction $\beta :B\to B\otimes A$ we have $\xi\in Z(B)\cap B^\beta$.

(ii) $\xi$ is a scalar if and only if $tr$ is the canonical trace.}

\lvs

{\bf PROOF.} By using a decomposition $B\simeq \bigoplus_{\gamma =1}^sM_{m_\gamma}$ we see that
$$\xi =\sum_\gamma m_\gamma^2\lambda_\gamma^{-1}Id_{M_{m_{\gamma}}}$$
(this is the last formula in the proof of proposition 2.1 in \cite{aut}). This proves both (ii) and the assertion $\xi\in Z(B)$ in (i). The remaining assertion $\xi\in B^\beta$ is clear from the fact that $\xi$ is a fixed vector of $u_\beta$ (cf. lemma 1.2 in \cite{aut}). 

\lvs

{\bf PROPOSITION 6.1.} {\it The index of a subfactor of the form $(B_0\otimes P)^{\beta_0\otimes\pi}\subset  (B_1\otimes P)^{\beta_1\otimes\pi}$ is equal to the index of the Markov inclusion $B_0\subset B_1$ and has to be an integer.}

\lvs

{\bf PROOF.} The first assertion is clear from the proof of theorem 5.1 (ii). Let $B_0=\bigoplus_iM_{n_i}$ and $B_1=\bigoplus_jM_{p_j}$ be decompositions of $B_0$ and $B_1$ as multimatrix algebras. Then the weights $\{\rho_i\}$ and $\{ \lambda_j\}$ of the traces of $B_0$ and $B_1$ are related by the formulae
$$\rho_i=\sum_j n_ip_j^{-1}m_{ij}\lambda_j$$
for any $i$, where $M=(m_{ij})$ is the inclusion matrix of $B_0\subset B_1$ (see \cite{ghj}). On the other hand we know from lemma 6.1 that the traces of $B_0$ and $B_1$ are their canonical traces, so their weights are $\rho_i=dim(B_0)^{-1}n_i^2$ and $\lambda_j=dim(B_1)^{-1}p_j^2$. We get that
$$dim(B_0)^{-1}n_i^2=\sum_j n_ip_j^{-1}m_{ij}dim(B_1)^{-1}p_j^2$$
i.e. that $\sum_j m_{ij}p_j=\gamma n_i$ for any $i$, with $\gamma =dim(B_1)/dim(B_0)$. This shows that $M^tp=\gamma n$, where $p$ and $n$ are the vectors $(p_j)$ and $(n_i)$. On the other hand from the definition of $M$ we have $Mn=p$. It follows that $(MM^t)p=\gamma p$. Since the index of $B_0\subset B_1$ is the unique Perron-Frobenius eigenvalue of $MM^t$, it is equal to $\gamma$. We know that $\gamma =dim(B_1)/dim(B_0)$ is a rational number, so we may write $\gamma =\alpha/\beta$ with $\alpha$ and $\beta$ positive integers having no common prime divisor. From $(MM^t)p=\alpha /\beta p$ we get succesively that $\beta$, then that $\beta^2$, then that $\beta^3$ etc. divides all $p_j$'s. Thus $\beta =1$. 

\vs

{\bf 7. APPENDIX: SEMIDUALITY OF MINIMAL COACTIONS}

\vs

A coaction $\pi :P\rightarrow P\otimes A$ is said to be minimal if $(P^\pi )^\prime\cap P={\bf C}$ and if it faithful in the following sense: ${\overline{sp}^w\,} \{ (\phi\otimes id)\pi (p)\mid \phi\in P_*,\, p\in P\} = A$ and semidual if each finite dimensional unitary corepresentation $u\in{\mathcal L} (H)\otimes A$ has a unitary eigenmatrix, i.e. a unitary $M\in {\mathcal L} (H)\otimes P$ such that $(id\otimes\pi )M=M_{12}u_{13}$. In this section we prove, following A. Wassermann, that the minimal coactions are semidual.

We use terminology from \S 2. Let $Irr(A)$ be a complete system of non-equivalent unitary irreducible corepresentations of $A$, each of them co-acting on some ${\bf C}^k$. Then the set $\{ dim(u)^{1/2}\,\widehat{u}_{ij}\}_{u\in Irr(A),i,j}$ is an orthonormal basis of the Hilbert space $l^2({\mathcal A} )$ (see \cite{wo}). For any $p\in {\mathcal P}$ we write $\pi (p)=\sum_{uij}p^u_{ij}\otimes u_{ij}$ and we use the following formulae:
$$\pi (p_{ij}^u)=\sum_k p_{kj}^u\otimes u_{ki},\,\,\,\,\, p=\sum_{ui}p^u_{ii},\,\,\,\,\, tr(p)=tr(p^1)$$

The formula for $\pi (p_{ij}^u)$ follows from the coassociativity of $\pi$. It shows that $\pi ({\mathcal P})\subset {\mathcal P}\otimes_{alg}{\mathcal A}$. By applying $1\otimes\varepsilon$ to the equality $(\pi_{\mid{\mathcal P}}\otimes id)\pi_{\mid{\mathcal P}}=(id\otimes\pi_{\mid{\mathcal P}} )\pi_{\mid{\mathcal P}}$ we get that $(\pi_{\mid{\mathcal P}}\otimes\varepsilon )\pi_{\mid{\mathcal P}} =\pi_{\mid{\mathcal P}}$, and since $\pi_{\mid{\mathcal P}}$ is injective it follows that $\pi_{\mid{\mathcal P}}$ is counital, i.e. that $(id\otimes\varepsilon )\pi_{\mid{\mathcal P}} =id$; this proves the formula for $p$. The formula for $tr(p)$ follows from the $\pi$-equivariance of the trace. 

\lvs

{\bf LEMMA 7.1.} {\it Let $\pi :P\to P\otimes A$ be a coaction. For any finite dimensional unitary corepresentation $u$ define a map $E_u:P\to P$ by 
$$E_u:p\mapsto dim(u)(id\otimes h)(\pi (p)(1\otimes \chi (u)^*))$$
Then $\{ E_u\}_{u\in Irr({\mathcal A} )}$ are orthogonal projections with respect to the trace of $P$, their images $P^u=Im(E_u)$ are in ${\mathcal P}=\pi^{-1}(P\otimes_{alg} {\mathcal A})$ and ${\mathcal P}$ decomposes as a direct sum $\oplus_u P^u$.}

\lvs

{\bf PROOF.} For any $p$ in the dense subalgebra ${\mathcal P}$ we have that
$$E_u(p)=dim(u)(id\otimes h)(\sum_{wij}(p^w_{ij}\otimes w_{ij})(1\otimes\chi (u)^*))=dim(u)\sum_{wijs} p^w_{ij} h(w_{ij}u_{ss}^*)=\sum_{i} p^u_{ii}$$
and together with $p=\sum_{ui}p^u_{ii}$ and $tr(p)=tr(p^1)$ this proves all the assertions. 

\lvs

{\bf LEMMA 7.2.} {\it Let $\alpha :P\to P\otimes A$ be a coaction and $u\in M_n\otimes{\mathcal A}$ be unitary corepresentation. Consider the unitary corepresentation
$$u^+:=(n\otimes 1)\oplus u=
\begin{pmatrix} 1&0\cr 0&u\end{pmatrix}
\in M_2(M_n\otimes {\mathcal A} )=M_2\otimes M_n\otimes {\mathcal A}$$
If the fixed point algebra $X_u=(M_2\otimes M_n\otimes P)^{\pi_{u^+}}$ is a factor then $u$ has a unitary eigenmatrix.}

\lvs

{\bf PROOF.} We have that $\pi_{u^+}$ leaves invariant $M_2\otimes M_n\otimes 1$ and the restriction is
$$\iota_{u^+}:x\otimes 1\mapsto \left( \begin{pmatrix} {1}&{0}\cr {0}&{u}\end{pmatrix}x\begin{pmatrix} {1}&{0}\cr {0}&{u^*}\end{pmatrix}
\right)_{124}$$
It follows that $X_u$ contains the two matrices $\begin{pmatrix}1&0\cr 0&0\end{pmatrix}$ and $\begin{pmatrix}0&0\cr 0&1\end{pmatrix}$, hence an element $K$ such that $KK^*=\begin{pmatrix}1&0\cr 0&0\end{pmatrix}$ and $K^*K=\begin{pmatrix}0&0\cr 0&1\end{pmatrix}$. Then $K$ has to be of the form  $\begin{pmatrix}0&m\cr 0&0\end{pmatrix}$ with $m\in M_n\otimes P$ unitary. The condition $K\in X$ is
$$\begin{pmatrix} {1}&{0}\cr {0}&{u_{13}}\end{pmatrix}
\begin{pmatrix} {0}&{(id\otimes\pi )m}\cr {0}&{0}\end{pmatrix}
\begin{pmatrix} {1}&{0}\cr {0}&{u^*_{13}}\end{pmatrix}
=\begin{pmatrix} {0}&{m\otimes 1}\cr {0}&{0}\end{pmatrix}$$
and this is equivalent to $(id\otimes\pi )m=m_{12}u_{13}$. 

\lvs

{\bf LEMMA 7.3.} {\it If $\pi :P\to P\otimes A$ is a minimal coaction and $u\in Irr(A)$ then $u$ has a unitary eigenmatrix if and only if $P^u \neq  \{ 0\}$.}

\lvs

{\bf PROOF.} The ``only if'' part follows from the fact that the entries of a $u$-eigenmatrix are in $P^u$. For the converse, by using lemma 7.2 it is enough to prove that the corresponding fixed point algebra $X_u$ is a factor. Let $x\in Z(X_u)$. We have $1\otimes
1\otimes P^\pi\subset X$ and from the irreducibility of $P^\pi\subset P$ we get that $x\in M_2\otimes M_n\otimes 1$. On the other hand we have 
$$X_u\cap M_2\otimes M_n\otimes 1=(M_2\otimes M_n)^{\iota_{u^+}}\otimes 1 =End(u^+)\otimes 1$$
Since $u$ is irreducible, it follows that $x$ is of the form $x=\begin{pmatrix} {y}&{0}\cr {0}&{\lambda I}\end{pmatrix}\otimes 1$ with $y\in M_n$ and $\lambda\in {\bf C}$. Let now $p\in P^u\neq \{ 0\}$ and write $\pi (p)=\sum_{ij}p_{ij}\otimes u_{ij}$. Then $\pi (p_{ij})=\sum_k p_{kj}\otimes u_{ki}$ for any $i,j$, i.e each column of $(p_{ij})_{ij}$ is a $u$-eigenvector. Choose such a non-zero column $l$ and let $m^i$ be the matrix having the $i$-th line equal to $l$ and zero elsewhere. Then $m_i$ is a $u$-eigenmatrix for any $i$, and this implies that  $\begin{pmatrix} {0}&{m^i}\cr {0}&{0}\end{pmatrix}\in X_u$ (cf. the end of the proof of lemma 7.2). The commutation relation of this matrix with $x$ is
$$\begin{pmatrix} {y}&{0}\cr {0}&{\lambda I}\end{pmatrix}
\begin{pmatrix} {0}&{m^i}\cr {0}&{0}\end{pmatrix}= 
\begin{pmatrix} {0}&{m^i}\cr {0}&{0}\end{pmatrix}
\begin{pmatrix} {y}&{0}\cr {0}&{\lambda I}\end{pmatrix}$$
and this gives $(y-\lambda I)m^i=0$, which by definition of $m^i$ shows that the $i$-th column of $y-\lambda I$ is zero. Thus $y-\lambda I= 0$, so $x=\lambda 1$. 

\lvs

{\bf THEOREM 7.1.} {\it The minimal coactions are semidual.}

\lvs

{\bf PROOF.} Let $K$ be the set of finite dimensional unitary corepresentations of ${\mathcal A}$ which have unitary eigenmatrices. Then:

-- $K$ is stable by tensor products: if $M$ (resp. $N$) is a unitary $u$- (resp. $w$-) eigenmatrix, then $M_{13}N_{23}$ is a unitary $u\otimes w$-eigenmatrix.

-- $K$ is stable by sums: if $M_i$ are unitary $u_i$-eigenmatrices then $diag(M_i)$ is a unitary eigenmatrix for $\oplus u_i$.

-- $K$ is stable by substractions: if $M$ is an eigenmatrix for $\oplus_{i=1}^{i=n} u_i$ then the first $dim(u_1)$ columns of $M$ are formed by elements of $P^{u_1}$, the next $dim(u_2)$ columns of $M$ are formed by elements of $P^{u_2}$ etc. Now if $M$ is unitary, it is in particular invertible, so all $P^{u_i}$'s are different from $\{ 0\}$ and we may conclude by using lemma 7.3.

-- $K$ is stable by complex conjugation: first, by the above results we may restrict attention to irreducible corepresentations. Now if $u\in Irr({\mathcal A} )$ has a nonzero eigenmatrix $M$ then $\overline{M}$ is an eigenmatrix for $\overline{u}$. By lemma 7.3 we get that $P^{\overline{u}}\neq\{ 0\}$, and we may conclude by using lemma 7.3.

These properties of $K$ and lemma 7.3 show that if ${\mathcal A}^\prime\subset{\mathcal A}$ is the $*$-subalgebra generated by the coefficients of the $u\in Irr({\mathcal A} )$ with $P^u\neq\{ 0\}$ then any corepresentation of ${\mathcal A}^\prime$ is in $K$. Together with the faithfulness assumption on $\pi$ this shows that ${\mathcal A}^\prime ={\mathcal A}$.

\baselineskip=12pt
\bigskip

\noindent Department of Mathematics, University of California, Berkeley, CA 94720, USA.

\noindent Institut de Math\'ematiques de Jussieu, 2 place Jussieu, 75005 Paris, France.

\noindent banica@math.berkeley.edu, banica@math.jussieu.fr.

\vs

\noindent AMS Classification: 46L37

\end{document}